\documentclass[11pt, reqno]{article}
\usepackage{amsfonts}
\usepackage{amsmath}
\usepackage{amsthm}
\usepackage{amssymb}
\def\b{\mathbb }

\def\phi{\varphi }

\swapnumbers
\theoremstyle{plain}
\newtheorem{theorem}{Theorem}[section]
\newtheorem{corollary}[theorem]{Corollary}
\newtheorem{lemma}[theorem]{Lemma}
\newtheorem{proposition}[theorem]{Proposition}
\theoremstyle{definition}
\newtheorem{definition}[theorem]{Definition}
\theoremstyle{remark}

\numberwithin{equation}{section}

\begin{document}
\title{Positivity of Dunkl's Intertwining Operator}
\author{Margit R\"osler\thanks{The main part of this paper was written
    while the  author 
held a Forschungsstipendium of the DFG at the University of Virginia, 
Charlottesville, VA,  U.S.A.} \\
 Zentrum Mathematik, Technische Universit\"at M\"unchen\\
 Arcisstr. 21, D-80290 M\"unchen, Germany\\
e-mail: roesler@mathematik.tu-muenchen.de
}

\date{}

\maketitle

\begin{abstract} For a  finite reflection group on $\b R^N,$ the  
  associated Dunkl operators  are parametrized first-order
  differential-difference operators  which
  generalize the usual  partial derivatives. They generate a commutative
  algebra which is -- under weak assumptions -- intertwined with the
  algebra of partial  differential operators by a unique linear and
  homogeneous isomorphism on polynomials.  
In this paper it is shown that for  non-negative parameter values, this
  intertwining operator is positivity-preserving on polynomials  and
  allows a positive integral representation on  certain algebras of
  analytic functions. This result in particular 
implies that the generalized exponential
  kernel of the Dunkl transform is positive-definite.  
\end{abstract}

\noindent
1991 AMS Subject Classification: Primary: 33C80; 
Secondary:  44A15, 33C50.

\smallskip

\noindent
Keywords and phrases: Dunkl operators, intertwining operator,
generalized exponential kernel.

\section{Introduction and results}

During the last years, the theory of Dunkl operators has found a wide area of 
applications in mathematics and mathematical physics. Besides their use in the 
study of multivariable orthogonality structures with certain reflection symmetries 
(see e.g. [D1-2], \cite{vD}, \cite{R}, \cite{X3}), these operators 
are for example closely related to certain 
representations of degenerate affine Hecke algebras (see \cite{Che}, \cite{Op} 
and, for some background, \cite{Ki}). Moreover, they have  
been succesfully involved in the description and solution of 
Calogero-Moser-Sutherland type quantum many body systems; among the wide
literature in this context, we refer to \cite{P}, \cite{LV} and \cite{BF}.

Let $G\subset O(N,\b R)$ be a finite reflection group on $\b R^N$. For
$\alpha\in \b R^N\setminus\{0\}$, we denote by $\sigma_\alpha $ the reflection in
the hyperplane orthogonal to $\alpha$, i.e.
\[\sigma_\alpha(x)\,=\, x -
2\,\frac{\langle\alpha,x\rangle}{|\alpha|^2}\,\alpha\,,\]
where $\langle.,.\rangle$ denotes the Euclidean
scalar product on $\b R^N$ and $|x|:=\sqrt{\langle x,x\rangle}$. (We
use the same notations for the standard Hermitian inner product and
norm on $\b C^N$.)
Let further $R$ be the  root system of $G$, normalized such that
$\langle\alpha,\alpha\rangle = 2$ for all $\alpha\in R$, and fix a
positive subsystem   $R_+$ of $R$. We recall from the general theory
of reflection groups (see e.g. \cite{Hu}) that the set of
 reflections in $G$ coincides with $\{\sigma_\alpha\,, \alpha\in
 R_+\}$, and that the orbits in $R$ under the natural action of $G$
 correspond to the conjugacy classes of reflections in $G$. A function
 $k: R\to \b C$ is called a multiplicity function on $R$, if it is
 $G$-invariant. We write $\text{Re}\,k\geq 0$ if
 $\text{Re}\,k(\alpha)\geq 0$ for all $\alpha\in R$ and $k\geq 0$ if
 $k(\alpha)\geq 0$ for all $\alpha\in R$. 

The Dunkl operators associated with $G$ are first-order 
differential-difference operators on $\b R^N$ 
which are parametrized by some 
multiplicity function $k$ on $R$. For $\xi\in \b R^N,$ the
corresponding Dunkl operator $T_\xi(k)$ is given by
\[ T_\xi(k)f(x) := \partial_\xi f(x) + \sum_{\alpha\in R_+} 
   k(\alpha)\,\langle\alpha, \xi\rangle\, 
    \frac{f(x) - f(\sigma_\alpha x)}{\langle\alpha, x\rangle}\,, 
    \quad\> f\in C^1(\b R^N);\]
here $\partial_\xi$ denotes the directional  derivative corresponding
    to $\xi$. As $k$ is
$G$-invariant, the above definition is independent of the choice of 
$R_+$. In case $k=0$, the $T_\xi(k)$ reduce to the 
corresponding directional
derivatives. The operators $T_\xi(k)$ were introduced and first studied by
Dunkl in  a series of papers ([D1-4]) in connection with a
generalization of the classical theory of spherical harmonics:
here the uniform surface measure on the $(N-1)$-dimensional unit
sphere is  modified by a weight function which is invariant under
the action of a given reflection group $G$ and associated with a
multiplicity function $k\geq 0$, namely
\begin{equation}\label{(1.1)} w_k(x) \,=\, \prod_{\alpha\in R_+}
    |\langle\alpha,x\rangle|^{2k(\alpha)}\,.\end{equation}
The most important properties of the operators $T_\xi(k)$ are as
follows: Let $\,\Pi^N = \b C[\b R^N]\,$ denote the algebra of
polynomial functions on $\b R^N$ and $\mathcal P_n^N$  ($n\in \b Z_+ =
\{0,1,2,\ldots\}$) the subspace of  homogeneous polynomials of (total)
degree $n$. 
Then  
\parskip=-2pt
\begin{enumerate}\itemsep=-1pt
\item[\rm{(1)}] The  $T_\xi(k)\,,\, \xi\in \b R^N,$ generate a
  commuting family of linear operators on $\Pi^N.$
\item[\rm{(2)}] Each $T_\xi(k)$ is homogeneous of degree $-1$ on
  $\Pi^N$, that is, $T_\xi(k)(p) \in \mathcal P_{n-1}^N$ for $p\in
  \mathcal P_n^N$. 
\item[\rm{(3)}] For all but a singular set of multiplicity functions,
  in particular for $k\geq 0$,  there exists a unique
  linear isomorphism $V_k$ of $\Pi^N$ such that 
  \[ V_k(\mathcal
  P_n^N)= \mathcal P_n^N\,,\>\> V_k\,\vert_{ \mathcal P_0^N}\, =\, id
  \quad\text{ and }\> 
  T_\xi(k) V_k\,=\, V_k\,\partial_\xi \quad\text{ for all }\,\xi\in \b R^N.\] 
\end{enumerate}
Properties (1) and (2) were  shown in \cite{Du1}, while the
existence of an intertwining operator according to (3) was first 
shown in \cite{Du2a} under the assumption $k\geq 0$.  An  abstract
and extended treatment of the above items is given in \cite{DJO}. 

The intertwining operator $V_k$ plays a central part in Dunkl's
theory and its applications. It is in particular involved in the definition of 
Dunkl's kernel  $K_G(x,y)$ (see below), which generalizes 
the usual exponential kernel $e^{\langle x, y\rangle}$ and arises as
the integral kernel of  the Dunkl transform (see \cite{Du3} and
\cite{deJ}). 
An explicit form of $V_k$ is known so far only
in very special cases: \parskip=-2pt
\begin{enumerate} \itemsep=-1pt
\item[\rm{1.}] The one-dimensional case, associated with the
reflection group $\b Z_2$ on $\b R$; here the
multiplicity function is given by a single parameter $k\geq 0$, and
the intertwining operator $V_k$ has the integral representation 
(see \cite{Du2}, Th. 5.1)
\begin{equation}\label{(1.9)}
 V_k p(x)\,=\, c_k \int_{-1}^1 p(xt)\,(1-t)^{k-1}(1+t)^k\,dt \quad
 \text{ with }\, c_k \,=\,
 \frac{\Gamma(k+1/2)}{\Gamma(1/2)\,\Gamma(k)}.
\end{equation} 
\item[\rm{2.}] The direct product case,
associated with the reflection group $\b Z_2^N$ on $\b R^N$; here a
closed form of the intertwining operator was determined in \cite{X1}.
\item[\rm{3.}] The case of the symmetric group $S_3$ on $\b R^3$,
  which has been studied in \cite{Du4}.
\end{enumerate} 

\medskip

In \cite{Du2}, the intertwining operator $V_k$ is, for $k\geq 0$,
extended to a
bounded linear operator on a suitably normed
algebra of  series of homogeneous polynomials on the unit ball.  To
allow a more convenient formulation of our
  statements, we introduce a slightly extended notation: For $r>0$ let
$\,K_r:=\{ x\in \b R^N: |x|\leq r\}$ denote the ball of radius $r$ and define
\begin{equation}\label{(1.10)}
 A_r := \big\{ f: K_r\to \b C , \, f=\sum_{n=0}^\infty f_n \quad\>\text{with
  }\> f_n\in \mathcal P_n^N \>\>\,\text{and}\>\,\> \|f\|_{A_r}:=
  \sum_{n=0}^\infty \|f_n\|_{\infty, K_r}\, <\infty \big\}  .
\end{equation}
It is easily checked that $A_r$ is a commutative Banach-*-algebra
(with complex conjugation as involution), see Section 4. Moreover, it 
follows from  Theorem 2.7 of \cite{Du2} that 
 $V_k$ extends to a continuous linear operator on  $A_r$  by
$\, V_k f := \sum_{n=0}^\infty V_k  f_n\,$ for $\, f=
\sum_{n=0}^\infty 
f_n \,\in A_r $. 
Up to now, it has been an open question whether for $k\geq 0$ the
intertwining operator $V_k$ is always positive, i.e. $V_kp \geq 0$ on
$\b R^N$ for each nonnegative polynomial $p\in \Pi^N$. More generally,
we may ask whether 
for every $x\in \b R^N$ with $|x|\leq r$, the functional 
$\, f\mapsto V_kf(x)$ is
positive on $A_r$.
 This property, which  was first conjectured in \cite{Du2} (in a
 slightly different setting), is
obvious in the above listed special cases 1 and 2 from the explicit
 representation of $V_k$;  in the $S_3$-case however,  the
integral representations derived in \cite{Du4} failed to infer
this result - at least for a large range of $k$. It is the aim of this
paper to prove the above conjecture for general reflection groups and
nonnegative multiplicity functions. Our first central result
establishes positivity of $V_k$ on polynomials: 

\begin{theorem} Assume that $k\geq 0$ and let $p\in \Pi^N$ with
  $p(x)\geq 0$ for all $x\in \b R^N$. 
Then also $V_k\,p(x)\geq 0$ for all $x\in \b R^N$.
\end{theorem}

More detailed information about $V_k$ is then obtained by its extension  to the
 algebras $A_r$. This 
 leads to the following theorem, which is the main result of this paper:

\begin{theorem} 
Assume that $k\geq 0$. Then for each $x\in \b R^N$ there exists a unique probability measure
$\mu_x$ on
the Borel-$\sigma$-algebra of $\b R^N$ such that 
\begin{equation}\label{(1.11)}
 V_kf(x)\,=\, \int_{\b R^N} f(\xi)\,d\mu_x(\xi)\quad \text{ for all}\>\, 
f\in A_{|x|}.
\end{equation} 
The representing measures $\mu_x$ are compactly
supported with
$\,{\rm supp}\,\mu_x \subseteq \, \{\xi\in \b R^N: |\xi|\leq |x|\}.\,$
Moreover, they satisfy
\begin{equation}\label{(1.12)}
\mu_{rx}(B)\,=\,\mu_x(r^{-1}B),\quad \mu_{gx}(B)\,=\,\mu_x(g^{-1}(B))
\end{equation} 
for each $r>0, \> g\in G$  and each Borel set $B\subseteq \b
R^N$.
\end{theorem}

\medskip

An important consequence of Theorem 1.2 concerns the generalized
exponential kernel $K_G,$ which is defined by 
\[ K_G(x,y) := V_k\bigl( e^{\langle\,.,\,y\,\rangle}\bigr)(x) \qquad (x,y \in
\b R^N),\]
see \cite{Du2}.
The function $K_G$ has a holomorphic extension to $\b C^N\times \b
C^N$ and is symmetric in its arguments. We also remark that by a
result of 
\cite{Op1}, the function $x\mapsto K_G(x,y)$ may  be characterized as
the unique analytic solution of the system $T_\xi(k)f =
\langle\xi,y\rangle\,f \>\>(\xi\in \b R^N)$ with $f(0)=1.$
 Theorem 1.2 implies that for fixed $y\in \b
R^N$ the kernel $ K_G(x,iy)$ is  positive-definite as a function of
$x$ on $\b R^N$, and the same holds for the "generalized Bessel
function" 
\[ J_G(x,iy) := \frac{1}{|G|} \sum_{g\in G} K_G(x,igy) \qquad(x,y\in \b
R^N).\]
As noted in \cite{Op1},  the kernel $J_G$ allows in some cases
(for Weyl groups $G$ and  certain discrete sets of multiplicity
functions)
an interpretation as the spherical function for some Euclidean
 symmetric space; in these cases positive-definiteness of $J_G$  is
 obvious. There are no similar
 interpretations known for the kernel $K_G$. Nevertheless, the conjecture that it should be positive-definite  has  been
confirmed by several of its properties (see \cite{deJ}), and  in particular by the fact that
$K_G(x,y)>0$ for all $x,y\in \b R^N$; this was proved in \cite{R}
in connection with the study of a generalized heat semigroup for Dunkl
operators. 
  
\medskip

The main parts of Theorem 1.2 are obtained by a standard
argumentation from Theorem 1.1. The proof of Theorem 1.1, however,  is much more involved.
Its crucial step  is a reduction from the 
$N$-dimensional to a one-dimensional problem,
 using semigroup techniques for 
linear operators on spaces of polynomials. The generators of the semigroups 
under consideration are certain differential-difference operators whose 
common decisive property is that they are ``degree-lowering''. This setting
is introduced in  Section 2, together with  a  Hille-Yosida type theorem 
which 
characterizes positivity of such semigroups by means of their generator. 
Theorem 1.1 is then  proved in 
Section 3. Section 4 is introduced   with  a short 
discussion of the algebras $A_r$ and their spectral properties, which is the 
basis for the subsequent  proof of Theorem 1.2. In the last section we
discuss  some implications of our results in the theory of Dunkl operators 
and related applications. 

\section{Semigroups generated by degree-lowering operators on polynomials}

We start with some 
general notations: Let $\,\Pi_+^N :=\{p\in \Pi^N\,: p(x)\geq 0$ for
all $x\in \b R^N\,\}\,$ denote the cone of nonnegative polynomials on 
$\b R^N$, and
$\Pi_n^N := \bigoplus_{k=0}^n \mathcal P_k^N \>\, (n\in \b Z_+)$ the space of
polynomials of (total) degree at most $n.$  The action of a subgroup $H
\subseteq O(N,\b R)$ on
$\Pi^N$ will always be the natural one, given by $\,hp(x):=
p(h^{-1}x)\> (h\in H, \>p\in \Pi^N)$. Finally, for a locally compact
Hausdorff space $X, \> M_b(X)$ is
the space of all regular bounded Borel measures on $X$ and
$M_b^+(X)$ the subspace of those which are non-negative.

\begin{definition}
A linear operator $A$ on $\Pi^N$ is called \parskip=-2pt
\begin{itemize}\itemsep=0pt
\item positive, if $Ap \in \Pi_+^N$ for each $p\in \Pi_+^N.$ 
\item degree-lowering, if 
$\,\displaystyle A\bigl(\Pi_n^N\bigr)\subseteq \Pi_{n-1}^N\,$ for all 
$n\in \b Z_+\,.$ 
\end{itemize}
\end{definition}

 \medskip

Important examples of degree-lowering operators are linear operators which 
are homogeneous of some degree $-n,\> n\geq 1$, on $\Pi^N$. This includes in 
particular  usual partial
derivatives and Dunkl operators, as well as products and linear combinations 
of those. 
If $A$ is  degree-lowering on $\Pi^N,$ then for every analytic 
function $f: \b R\to \b C$ with power series $f(x)= \sum_{k=0}^\infty 
c_k x^k$, there is a linear operator $f(A)$ on $\Pi^N$ defined by 
the terminating series 
\[ f(A) p\,(x):= \sum_{k=0}^\infty c_k A^k p\,(x). \]
Notice that $\,f(A)(\Pi_n^N)\subseteq \Pi_n^N$ for 
each $n\in \b Z_+\,.$ This yields a natural restriction of 
$f(A)$ to a linear operator on the 
finite-dimensional vector space $\Pi_n^N$. In particular, the well-known 
product and exponential formulas for linear operators on finite-dimensional 
vector spaces (see, e.g. \S 4.7 of \cite{K}) imply corresponding exponential 
formulas 
for degree-lowering operators on  $\Pi^N$, where the topology may be 
choosen to be the one of pointwise convergence.  We note two results of 
this type, which will be used later on:

\begin{lemma}
Let $A$ and $B$ be degree-lowering linear operators on $\Pi^N$. Then for
 all $p\in \Pi^N$ and $x\in \b R^N$, \parskip=-1pt
\begin{enumerate} 
\item[\rm{(i)}] $\displaystyle e^A \,p\,(x)\,=\, \lim_{n\to\infty} 
  \Bigl(I-\frac{A}{n}\Bigr)^{-n} p\,(x).$
\item[\rm{(ii)}] $\displaystyle e^{A+B} p\,(x)\,=\, \lim_{n\to\infty}
 \bigl(e^{A/n} e^{B/n}\bigr)^n p\,(x).\quad$ (Trotter product formula).
\end{enumerate}
\end{lemma}

Each degree-lowering  operator $A$ on $\Pi^N$ generates a 
semigroup $\,(e^{tA})_{t\geq 0}\,$ of linear operators on $\Pi^N$ and, 
in fact, on each of the $\Pi_n^N.$ 
The generator $A$ is uniquely determined from the semigroup by 
\[ Ap\,(x)\,=\, \lim_{t\downarrow 0}\, t^{-1}\bigl(e^{tA}-I\bigr)p\,
  (x)\quad \text{for all }\, p\in \Pi^N.\]

The following key-result characterizes positive semigroups generated 
by degree-lowering operators; it is an adaption of a well-known Hille-Yosida 
type 
characterization theorem for Feller-Markov semigroups on $C(K),\> K$ a compact
 Hausdorff space (see, e.g.  \S II.4 of \cite{GS}):

\begin{theorem}
Let $A$ be a degree-lowering linear operator on $\Pi^N$. Then these are
 equivalent: \parskip=-2pt
\begin{enumerate}\itemsep=1pt
\item[\rm{(1)}] $\, e^{tA}$ is positive on $\Pi^N$ for all $t\geq 0$.
\item[\rm{(2)}] $A$ satisfies the ``positive minimum principle''
\parskip=-1pt
\begin{center}
 {\rm (M)} \hspace{2pt} For every $p\in \Pi_+^N$ and $x_0\in \b R^N,$  
\hspace{2pt} 
  $p(x_0)=0\>$ implies $\, Ap\,(x_0)\geq 0.$ 
\end{center}
\end{enumerate} 
\end{theorem}

\begin{proof}
(1) $\Rightarrow$ (2): Let $p\in \Pi_+^N$ with $p(x_0)=0.$ Then 
\[ Ap\,(x_0)\,=\, \lim_{t\downarrow 0} \frac{e^{tA}p\,(x_0) - p(x_0)}{t}\,\,=
\,\, 
   \lim_{t\downarrow 0}\frac{1}{t}\, e^{tA} p\,(x_0)\,\geq 0.\]

\noindent
(2) $\Rightarrow$ (1): Notice first that for each $\lambda \not= 0$, the 
operator $\lambda I -A$ is bijective on $\Pi^N.$ In fact, $\lambda I -A$ is 
injective on $\Pi^N$, because otherwise there would exist some $p\in 
\Pi^N,\>p\not=0$, with $\,Ap=\lambda p,$ in contradiction to the 
degree-lowering character of $A$. As $(\lambda I -A)(\Pi_n^N)\subseteq 
\Pi_n^N$, this already proves bijectivity of $\lambda I -A$ on each $\Pi_n^N$, 
hence on $\Pi^N$ as well. We next claim that for every $\lambda >0$ 
the resolvent operator $\,R(\lambda;A):= (\lambda I -A)^{-1}$  is 
positive on 
$\Pi^N$. For this, let $p\in \Pi^N_+$ and $q:=R(\lambda;A)\,p.$ If $p$
is constant, 
then $q=\frac{1}{\lambda}p \geq 0$. We may therefore
restrict to the case that the total degree $n$ of $p$ (which must be even)
is greater than $0$.
Suppose first that $p(x)\geq c|x|^n\,$ for all
$x\in \b R^N$, with some constant $c>0$. As $A$ lowers the degree, we
may write 
$\,q= \frac{1}{\lambda} p + r\,$ with a polynomial $r$ of
total degree less than $n$. Hence $\,\lim_{|x|\to\infty}
q(x)=\infty$, which shows 
that $q$ takes an absolute minimum, let us say in $x_0\in \b R^N$. Put 
$\,\widetilde q(x):= q(x)-q(x_0).$ Then $\widetilde q\in \Pi_+^N$ with 
$\widetilde q(x_0)=0,$ and property (M) assures that 
$\,Aq(x_0)= A\widetilde q(x_0) \geq 0.$ For $\lambda>0$ and $x\in \b R^N$ we 
therefore obtain
\[ \lambda q(x)\,\geq\, \lambda q(x_0)\,=\, (\lambda I -A)q(x_0) + Aq(x_0)\,
\geq\, p(x_0)\,\geq 0.\]
If $p\in \Pi_+^N$ is arbitrary, then consider the polynomials
$\,p_\epsilon(x):= p(x) +\epsilon\,|x|^n\,$ for $\epsilon>0$, where $n$ is the
degree of $p$. As $A$ is degree-lowering, and by the above
result, we obtain
\[ R(\lambda;A)\,p\,(x) \,=\, \lim_{\epsilon\to 0}
R(\lambda;A)\,p_\epsilon(x)\,\geq 0 \quad\text{for all}\>\> x\in \b R^N.\]
This proves the stated positivity of $R(\lambda;A)$ for $\lambda>0.$
\noindent
Now let $p\in \Pi_+^N$ and $t >0$.  Then according to Lemma 2.2.(i),
\[e^{tA} p\,(x)\,=\, \lim_{n\to \infty} \Bigl( I- \frac{tA}{n}\Bigr)^{-n} 
p\,(x)
 \,=\,\lim_{n\to\infty} \Bigl(\frac{n}{t} \,R\bigl(\frac{n}{t}; 
A\bigr)\Bigr)^n p\,(x)\,\geq 0\]
for all $x\in \b R^N$. This finishes the proof.
\end{proof}

\section{Positivity of $V_k$ on polynomials}

This section is devoted to the proof of  Theorem 1.1. The outline of
this proof is as follows: 
In a first step, we consider the 
(one-dimensional) differential-difference operators 
\[ \Lambda_s :=\, e^{-sD^2} \delta \,e^{sD^2}, \quad s\geq 0 \]
on $\Pi^1.$ Here $D$ denotes the usual first derivative, i.e. 
$Dp(x) = p^\prime (x)$ for $x\in \b R$, and $\delta$ is the linear operator on 
$\Pi^1$
 given by
\begin{equation}\label{(3.30)} 
\delta p(x)\,=\, \frac{p^\prime(x)}{x} \,-\, \frac{p(x)-p(-x)}{2x^2}\,
=\,\, \frac{1}{2}\int_{-1}^1 \bigl(D^2p\bigr)(tx)(1+t)\,dt.
\end{equation}
This operator is related to  the Dunkl operator $T(k)$ associated with 
the reflection group $\b Z_2$ on $\b R$ and the multiplicity parameter 
$k \geq 0$  by  
\[ T(k)^2\,=\, D^2 + 2k\delta.\]
As both $D^2$ and $\delta$ are homogeneous of degree $-2$  on $\Pi^1$, 
the operators
$\Lambda_s$ are well-defined and  degree-lowering on $\Pi^1$. We shall prove
that they have the following decisive property:

\begin{proposition}
The operators $\Lambda_s\,, \> s\geq 0$, satisfy the positive minimum principle 
{\rm (M)} on 
$\Pi^1.$  
\end{proposition}

We next turn to the general $N$-dimensional setting: Here
$G$ is an arbitrary finite reflection group on $\b R^N$ with multiplicity
function $k\geq 0$. We consider the generalized Laplacian associated
with $G$ and $k$, which is defined by 
\[\Delta_k := \sum_{i=1}^N T_{\xi_i}^2\,\]
 with an arbitrary orthonormal basis
$(\xi_1,\ldots,\xi_n)$ of $\b R^N$ (see \cite{Du1}). 
It is homogeneous of degree $-2$
on $\Pi^N$ and (with our convention $\langle\alpha,\alpha\rangle = 2$ for all
$\alpha\in R_+$) given explicitly by 
\begin{equation}\label{(3.20)}
 \Delta_k \,=\, \Delta  + 2\sum_{\alpha\in R_+}
k(\alpha)\delta_\alpha \quad\,\text{ with }\quad \delta_\alpha f(x)\,=\,
\frac{\langle\nabla f(x),\alpha\rangle}{\langle\alpha ,x\rangle} - 
\frac{f(x)-f(\sigma_\alpha x)}{\langle\alpha , x\rangle^2};
\end{equation}
 here $\Delta$ and $\nabla$ denote the usual Laplacian and gradient 
respectively. Theorem 2.3 is  the key to infer from the
one-dimensional setting of Proposition 3.1 to a general
multi-variable extension:

\begin{proposition}
Let $\, L_k := \Delta_k - \Delta. $ 
Then for $k\geq 0,$ the operators 
\[ e^{-s\Delta} e^{tL_k} e^{s\Delta}\,\quad (s,t\geq 0)\]
are positive on $\Pi^N.$ 
\end{proposition}

The statement of Theorem 1.1 will then finally be reduced to the following 
consequence of Proposition 3.2:

\begin{corollary}
The operator $\,\displaystyle e^{-\Delta/2} e^{\Delta_k/2}\,$ is  
positive on $\Pi^N$. 
\end{corollary}

\begin{proof} Applying Trotter's product formula of Lemma 2.2, we obtain 
\begin{align}
 e^{-\Delta/2} e^{\Delta_k/2}\, p\,(x)\,=\,&\, e^{-\Delta/2} e^{\Delta/2 + 
   L_k/2}\, p\,(x)\,= \,
  \lim_{n\to\infty}\,e^{-\Delta/2}\Bigl( e^{\Delta/2n}\, 
   e^{L_k/2n}\Bigr)^n p\,(x)\, \notag\\
 =\,& \lim_{n\to\infty} \prod_{j=1}^n \Bigl(e^{-(1-j/n)\cdot\Delta/2}\,
    e^{L_k/2n}\,e^{(1-j/n)\cdot\Delta/2}\Bigr) p\,(x) \quad\>\> 
   (p\in \Pi^N, \> x\in \b R^N).\notag 
 \end{align}
By Proposition 3.2, each of the $n$ factors in the above product is a 
positive operator on $\Pi^N$. Hence $\,e^{-\Delta/2} e^{\Delta_k/2}\,$ is also
 positive on $\Pi^N$. 
\end{proof}
 
We now turn to the proof of Proposition 3.1. We start with two elementary 
auxiliary results:

\begin{lemma}
For each $p\in \Pi^1$ and $c\in \b R$, 
\[ e^{cD^2}(xp(x))\,=\, x\,e^{cD^2}\!p\,(x) + 2c\,e^{cD^2} p^\prime (x).\]
\end{lemma}

\begin{proof}
Power series expansion of $e^{cD^2}$ yields
\begin{align}
e^{cD^2}(xp(x))\,=\,& \sum_{n=0}^\infty \frac{c^n}{n!}\, D^{2n}(xp(x))\,=\, 
xp(x) + \sum_{n=1}^\infty \frac{c^n}{n!} \bigl(x\, D^{2n}p\,(x) + 
2n D^{2n-1}p(x)\bigr)\notag\\
=\, & x\, e^{cD^2} p\,(x) +\,2c\sum_{n=1}^\infty \frac{c^{n-1}}{(n-1)!}\, 
D^{2n-1}p\,(x)\,=\, 
x\, e^{cD^2} p\,(x) +\,2c\, e^{cD^2}p^\prime(x). \notag
\end{align}
\end{proof}

\begin{lemma}
Let $p\in \Pi^1_{2n+1}\,,\>\, n\in \b Z_+\,,$ be an odd polynomial. Then 
the differential equation
\begin{equation}\label{(3.1)} 
 c\,y^\prime -xy\,=\, p \qquad  (c>0)
\end{equation}
has exactly one polynomial solution (which belongs to $\Pi^1_{2n}$), namely
\[ y_p(x)\,=\, \frac{1}{c}\, e^{x^2/2c}\int_{-\infty}^x e^{-t^2/2c}\, 
p(t)\,dt.\] 
\end{lemma}

\begin{proof}
The general solution of \eqref{(3.1)} is
\[ y(x)\,=\, a\,e^{x^2/2c}\,+\,\frac{1}{c}\,e^{x^2/2c}\int_{-\infty}^x 
e^{-t^2/2c}\,p(t)\,dt, \quad a\in \b R.\]
It therefore remains to prove that 
\begin{equation}\label{(3.2)}
x\mapsto e^{x^2/2c} \int_{-\infty}^x e^{-t^2/2c}\, p(t)\,dt 
\end{equation}
is a polynomial. We use induction by $n$: 
For $n=0$, the statement is obvious. For $n\geq 1,$ 
write $p(x)= -c^{-1}xr(x)$ with $r\in \Pi^N_{2n}.$ Partial 
integration then yields
 \[\int_{-\infty}^x e^{-t^2/2c}\,p(t)\,dt\,=\, -\frac{1}{c}\int_{-\infty}^x t
e^{-t^2/2c} \,r(t)\,dt\,=\,
e^{-x^2/2c}\, r(x)\,-\, \int_{-\infty}^x e^{-t^2/2c}\, r^\prime(t)\,dt.
\]
By our induction hypothesis, this equals
$\, e^{-x^2/2c} (r(x)-\widetilde r(x))\,$
with some polynomial $\widetilde r\in \Pi^N_{2n-2}.$  This finishes the proof.
\end{proof}

\begin{proof}[{\bf Proof of Proposition 3.1}] \, 1. The case $s=0$ is easy and 
may be treated separately:
Let $p\in \Pi^1_+$ with $p(x_0)=0.$ Then $p^\prime (x_0)=0$ and 
$p^{\prime\prime}(x_0)\geq 0.$ Thus if $x_0\not=0$, then $\,\delta p\,(x_0) = 
p(-x_0)/(2x_0^2) \geq 0.$ In case $x_0=0$, it is seen from the  integral 
representation 
\eqref{(3.30)}  that $\delta p\,(0) = p^{\prime\prime}(0) \geq 0.$ 

\noindent
From now on, we may therefore assume that $s>0$.    

\noindent
2. We first derive an explicit 
representation of the operator $\Lambda_s\>\> (s>0),$ which allows to check 
property (M) easily: We claim that
\begin{align}\label{(3.3)} 
\Lambda_s p\,(x)\,=\,& -\frac{1}{2s} p(x)\,-\, \frac{1}{4s^2}\, e^{x^2/4s}\,
\Bigl(\int_{-\infty}^x g_{p,\,x}(t)\,dt\,-\,\int_{-x}^\infty 
g_{p,\,x}(t)\,dt\Bigr) \quad\text{ for }\,p\in \Pi^1,\notag\\
 &\text{with}\quad g_{p,\,x}(t)\,=\, e^{-t^2/4s}\,(t+x)\,p(t). 
\end{align}
This may of course be verified by a (tedious) direct computation of 
$\Lambda_s(x^k),
 \>\, k\in \b Z_+$, and an explicit evaluation of the corresponding integrals 
on the right side by series expansions of the involved exponentials. We 
prefer, however, to give a more instructive proof:  

Note first that the operators $D^2$ and $\delta$ map even polynomials to even 
ones and odd polynomials to odd ones again, and that
\[ \delta p\,(x)\,=\,\begin{cases}
        \displaystyle{\frac{1}{x} p^\prime(x)}& \text{if $p$ is even},\\
        \displaystyle{\Bigl(\frac{1}{x}\,p(x)\Bigr)^\prime}& \text{if $p$ is 
                                                    odd}.
          \end{cases}
\]
Now fix $s>0$ and suppose that $p\in \Pi^1$ is even. Then the polynomials 
$\,e^{sD^2}\!p$ and $\,q:= \Lambda_sp$ are also even, and we obtain the following 
equivalences:
\[ q=\Lambda_sp\, \Longleftrightarrow\, \delta\bigl(e^{sD^2} \!p\bigr)\,=\, 
e^{sD^2}\!q\,
\Longleftrightarrow\,p^\prime(x)\,=\,e^{-sD^2}\bigl( x\,e^{sD^2}\!q\bigr)(x).\]
By use  of Lemma 3.4, this becomes
\begin{equation}
 p^\prime(x)\,=\, xq(x)\,-\, 2s q^\prime(x),
\end{equation}
which is a differential equation of type \eqref{(3.1)} for $q$.
Lemma 3.5, together with a further partial integration, now implies that 
\begin{align}\label{(3.4)}
 \Lambda_sp\,(x)\,=\,& -\frac{1}{2s}\,e^{x^2/4s} \int_{-\infty}^x 
e^{-t^2/4s}\, 
p^\prime(t)\,dt \notag\\
=\, & -\frac{1}{2s}\,p(x)\,-\, \frac{1}{4s^2}\,e^{x^2/4s}\int_{-\infty}^x 
e^{-t^2/4s}\, t\, p(t)\,dt \qquad\text{($p$ even).}
\end{align}
In a similar way, we calculate $\,q=\Lambda_sp$ for odd $p\in \Pi^1$:
In this 
case, 
$\,e^{sD^2}\!p\,$ and $q= \Lambda_s p$ are odd as well, and we have
the 
equivalence 
\[ q= \Lambda_sp \, \Longleftrightarrow\, \frac{d}{dx} 
   \Bigl(\frac{1}{x}\,e^{sD^2} p\,(x)\Bigr)\,=\, e^{sD^2}\! q\,(x).\]
Hence there exists a constant $c_1\in \b R$ such that 
\[ e^{sD^2}\! p\,(x)\,=\, x(c_1 + h(x)), \quad\text{ with }\>\, h(x)\,=\, 
   \int_0^x e^{sD^2}\!q\,(t)\,dt.\]
Applying Lemma 3.4 again, we obtain
\begin{equation}\label{(3.5)}
 p(x)\,=\, c_1 e^{-sD^2}(x) + x e^{-sD^2}h\,(x) - 2s\,e^{-sD^2}h^\prime
   (x)\,=\, c_1x + x\,e^{-sD^2}\! h\,(x) - 2sq(x). 
\end{equation} 
In order to determine $\,e^{-sD^2}\!h$, note that
\[\frac{d}{dx} \bigl(e^{-sD^2}\!h\,(x)\bigr)\,=\, e^{-sD^2} h^\prime (x)\,=\,
q(x).\]
Consequently, there exists a constant $c_2\in \b R$ such that
\begin{equation}\label{(3.6)}
 e^{-sD^2} h\,(x)\,=\, c_2 + \int_0^x q(t)\,dt.
\end{equation}
Now write $\,p(x)= x\,P(x)\,$ and $\,q(x) = x\, Q(x)\,$ with even $P,\,Q\in 
\Pi^1$. Then by \eqref{(3.5)} and \eqref{(3.6)},
\[ P(x)\,=\, c_1 + c_2 + \int_0^x t\,Q(t)\,dt\,- 2s\,Q(x),\]
and therefore
\[ P^\prime(x)\,=\, x\,Q(x) - 2s\,Q^\prime(x).\]
This is exactly the same differential equation as we had  in the 
even case before, and the transfer of \eqref{(3.4)} gives
\begin{equation}\label{(3.7)} 
  \Lambda_sp\,(x)\,=\, -\frac{1}{2s}\, p(x)\,-\, \frac{1}{4s^2}\, e^{x^2/4s}\,
  x\int_{-\infty}^x e^{-t^2/4s} p(t)\,dt \qquad \text{($p$ odd).}
\end{equation}
If finally $p\in \Pi^1$ is arbitrary, then write $p=p_e + p_o\,$ with even 
part $\,p_e(x) = (p(x)+ p(-x))/2\,$ and odd part $\,p_o(x)= (p(x)-p(-x))/x\,.$
The combination of \eqref{(3.4)} for $p_e$ with \eqref{(3.7)} for $p_o$ then 
leads to
\[ \Lambda_sp\,(x)\,=\, -\frac{1}{2s}\,p(x)\,-\, \frac{1}{4s^2}\, e^{x^2/4s}
\int_{-\infty}^x e^{-t^2/4s} \Bigl( \frac{t+x}{2}\,p(t)\,+\,\frac{t-x}{2}\, 
p(-t)\Bigr)dt,\]
and an easy reformulation yields the stated representation \eqref{(3.3)}.

\noindent
3. In order to prove that $\Lambda_s$ satisfies the positive minimum 
principle (M),
 define
\[ F_p(x):= \int_{-\infty}^x g_{p,\,x}(t)\,dt\,-\, \int_{-x}^\infty 
g_{p,\,x}(t)\,
dt,\quad\text{for }\,p\in \Pi^1 \>\>\text{and }\,x\in \b R\,.\]
Now let $p\in \Pi^1_+$ with $p(x_0)=0.$ Then in view of  \eqref{(3.3)}, 
\[ \Lambda_sp\,(x_0)\,=\, -\frac{1}{4s^2} \, e^{x_0^2/4s}\, F_p(x_0),\]
and it remains to check that $F_p(x_0)\leq 0.$ For this, we rewrite $F_p$ as
\[ F_p(x)\,=\, \int_{-\infty}^{-|x|} g_{p,\,x}(t)\,dt\,-\, \int_{|x|}^\infty 
g_{p,\,x}(t)\,dt.\]
As $p$ is nonnegative, the sign of $\,g_{p,\,x}(t)$ coincides with the sign of 
$(x+t)$ for all $x,t\in \b R.$ This shows that in fact, $\,F_p(x)\leq 0$ for 
all $x\in \b R$, which completes the proof. 
\end{proof}

\medskip

\begin{proof}[{\bf Proof of Proposition 3.2}] For fixed $s\geq 0,$ the 
operators $\,\bigl(e^{-s\Delta} e^{t\,L_k} e^{s\Delta}\bigr)_{t\geq 0} \,$ 
form a semigroup on $\Pi^N$ with generator 
$\,\displaystyle e^{-s\Delta} L_k \,e^{s\Delta}.\,$
According to Theorem 2.3, it therefore suffices to prove that this generator 
satisfies the  positive minimum principle (M) on $\Pi^N$. With the
notation of \eqref{(3.20)}, we have 
\[L_k = \,2\sum_{\alpha\in R_+} k(\alpha)\, \delta_\alpha
\quad\text{ with }\,\>
k(\alpha)\geq 0 \quad\>\text{for all }\,\alpha\in R_+\,. \]
It is therefore enough to make sure that each of 
the operators 
\[ \rho^s_{\alpha} := \, e^{-s\Delta} \,\delta_\alpha\, e^{s\Delta} 
\quad \> (\alpha\in R_+)\]
satisfies (M). (Here the assumption $k\geq 0$ is crucial!) Now fix $\alpha\in R_+$. An easy calculation shows that 
$\delta_\alpha$ and hence also $\rho^s_{\alpha}\,$  is 
rotation-equivariant, i.e.
\[ g\circ \rho^s_{\alpha} \circ g^{-1}\,=\, \rho^s_{g(\alpha)} \quad \>
\text{for }\, g\in {SO}(N,\b R).\]
We may therefore assume that $\,\alpha = 
\sqrt 2 e_1 = (\sqrt 2, 0,\ldots, 0).$ As  $\,\delta_{\sqrt 2 e_1}\,$ 
obviously commutes with each of the partial derivatives $\partial_2\,,\ldots,
 \partial_N\,$ on $\b R^N$, we obtain
\[ \rho^s_{\sqrt 2 e_1}\,=\, e^{-s\partial_1^2}\,
   \delta_{\sqrt 2 e_1}\,e^{s\partial_1^2}\,.\]
But this operator acts on the first variable only, namely via $\Lambda_s$:
\[ \rho^s_{\sqrt 2 e_1}\, p\,(x_1\,,\ldots, x_N)\,=\, 
\Lambda_s \,p_{x_2\,\ldots,\,x_N}\,(x_1), \>\> \text{where }\>\, 
   p_{x_2\,,\ldots,\,x_N}\,(x_1):= p(x_1,x_2,\ldots,x_N) \>\>\text {for } 
    p\in \Pi^N.\]
The assertion now follows from Proposition 3.1.
\end{proof}

In order to complete the proof of Theorem 1.1, we employ the following
bilinear form on $\Pi^N$ associated with $G$ and $k$, which was
introduced in \cite{Du2} (for a further discussion, see also
\cite{DJO}):
\[ [p,q]_k \,:=\, \bigl(p(T_k)\,q\bigr)(0) \quad\text{for } p,q\in
\Pi^N\,;\]
here $p(T_k)$ is the differential-difference operator which is obtained
from $p(x)$ by replacing each $x_i$ by the corresponding Dunkl
operator $T_{e_i}(k)$. The case $k=0$ will be distinguished by the notation
$p(\partial)$. Notice that $[p,q]_k = 0\,$ for 
$p\in \mathcal P_n^N$ and $q\in \mathcal
P_m^N$ with $n\not=m$.  It was shown in \cite{Du2} that for $k\geq 0$
 and  for all $p,q\in \Pi^N$,
\begin{equation}\label{(3.32)}
 \big[p,q\big]_k \,=\, c_k\int_{\b R^N} e^{-\Delta_k/2} p\,(x)\, 
   e^{-\Delta_k/2} q\,(x)\,e^{-|x|^2/2}\,w_k(x) dx,
\end{equation}
where $w_k$ is the weight function defined in \eqref{(1.1)} and $\displaystyle 
c_k := \Bigl(\int_{\b R^N}  e^{-|x|^2/2}\,w_k(x) dx\Bigr)^{-1}\,.$ We
remark that \eqref{(3.32)} can also be proved in a completely
independent way by using certain biorthogonal polynomial systems
(Appell characters and cocharacters) in $L^2(\b R^N,
e^{-|x|^2/2}w_k(x)dx)$; see \cite{RV1}. 
Another useful identity for $[.\,,.]_k$ is
\begin{equation}\label{(3.31)} 
[V_k\, p,q]_k\,=\, [p,q]_0 \quad\text{ for all }\,p,q\in
\Pi^N.
\end{equation}
In fact, for $p,q\in \mathcal P_n^N$ with  $n\in \b Z_+$, one obtains
\[
[V_k\, p,q]_k\,=\, [q,V_k\, p]_k\,=\,q(T_k)(V_k\, p)\,=\, 
   V_k(q(\partial)p)\,=\, q(\partial)(p)\,=\,
[p,q]_0\,; \]
here the characterizing properties of $V_k$ and the fact that
   $q(\partial)p$ is a constant have been used. 
For general $p,q\in \Pi^N$, \eqref{(3.31)} then follows from the
orthogonality of the spaces $\mathcal P_n^N\,,\> n\in \b Z_+\,,$ with
respect to both scalar products.  

Finally, we shall need the following positivity criterion for polynomials:

\begin{lemma}
Let $\alpha > 0$ and suppose that $\,h\in C_b(\b R^N)$ satisfies
\begin{equation}\label{(3.8)}
 \int_{\b R^N} h(x)\,p(x)\,  e^{-\alpha|x|^2}\,w_k(x)dx \, \geq 0 \quad\>\> 
   \text{for all }\, p\in \Pi^N_+\,.
\end{equation}
Then $\,h(x)\geq 0\,$ for all $x\in \b R^N$. 
\end{lemma}

\begin{proof} For abbreviation, put 
\[ dm_k(x):=  e^{-\alpha |x|^2}\,w_k(x)dx \,\in M_b^+(\b R^N)\] 

\noindent
1. We shall use that $\Pi^N$ is dense in $L^2(\b R^N, dm_k)$. This is 
proved (with
$\alpha= 1/2$) in Theorem 2.5 of \cite{Du3} by referring to a
well-known Theorem 
of Hamburger for one-dimensional distributions, but it can also be seen 
directly as follows: Suppose that  $\Pi^N$ is not dense in
$L^2(\b R^N, dm_k).$
Then there exists some $f\in L^2(\b R^N, dm_k), \> f\not= 0$, with 
$\,\int_{\b R^N} fp\,d m_k\,= 0\,$ for all $\,p\in \Pi^N$. Now consider
the measure $\nu:= fm_k\,\in M_b(\b R^N)$ and its (classical) 
Fourier-Stieltjes transform 
\[ \widehat\nu(\lambda)\,=\, \int_{\b R^N} e^{-i\,\langle\lambda,x\rangle}\, 
d\nu(x)\,=\, \int_{\b R^N} f(x)\, e^{-i\,\langle \lambda, x\rangle}\, 
dm_k(x).\]
As $\,x\mapsto e^{|\lambda||x|}\,$ belongs to $L^2(\b R^N, dm_k)$ for all 
$\lambda\in \b R^N$, the dominated convergence theorem yields  
\[ \widehat \nu(\lambda)\,=\, \sum_{n=0}^\infty \frac{(-i)^n}{n!}
  \int_{\b R^N} f(x)\,\langle\lambda,x\rangle^n \,dm_k(x)\,=\,0. \]
By injectivity of the Fourier-Stieltjes transform on $M_b(\b R^N)$, it 
follows that $\nu=0$ and hence $f=0$ a.e., a contradiction. 

\noindent
2. Now assume that $h\in C_b(\b R^N)$ satisfies \eqref{(3.8)}. In
order to
prove  $h\geq 0$, it suffices to check that  
\begin{equation}\label{(3.80)}
\int_{\b R^N} fh\,dm_k\,\geq 0 \quad\>\>\, \text{for all }\,\, 
f\in C_c^+(\b R^N).
\end{equation}
For this, let $f\in C_c^+(\b R^N)$ and $\epsilon>0$. By density of $\Pi^N$ in 
$L^2(\b R^N, dm_k)$ there exists some $p = p_\epsilon \in \Pi^N$ 
with $\,\displaystyle \|\sqrt f - p\,\|_{2,m_k} <\epsilon$.  With $\,M:=
\|h\|_{\infty,\b R^N}$ it follows that
\begin{align}
 \Big|\int_{\b R^N} fh\,dm_k\,-\,& \int_{\b R^N} p^2 h\,dm_k\Big|\,\leq\, 
M\int_{\b R^N} |f- p^2|\, dm_k\,\notag\\
\leq\, & \,M\cdot \|\sqrt f - p\,\|_{2,m_k} \,\,\|\sqrt f +
p\,\|_{2,m_k} \,
\leq\, M\epsilon \cdot\bigl(2\,\|\sqrt f\,\|_{2,m_k} + \epsilon\bigr),\notag
\end{align}
 which tends to $0$ with $\epsilon\to 0$. This proves \eqref{(3.80)}
 and yields the assertion.
\end{proof}

The proof of Theorem 1.1 is now easily accomplished:
 
\begin{proof}[{\bf Proof of Theorem 1.1}]  
Combining formulas \eqref{(3.31)} and \eqref{(3.32)}, we obtain for all $p,q\in
\Pi^N$ the identity
\[
c_k\int_{\b R^N} e^{-\Delta_k/2}(V_k\, p)(x)\, e^{-\Delta_k/2} q(x)\,
e^{-|x|^2/2}\,w_k(x) dx\,=\, c_0\int_{\b R^N} e^{-\Delta/2}
p(x)\,e^{-\Delta/2} q(x)\,e^{-|x|^2/2} dx.
\]
As $\displaystyle\, e^{-\Delta_k/2}(V_k\, p)\,=\,
V_k\bigl(e^{-\Delta/2}p\bigr)\,$, and as we may also replace $p$ by
$e^{\Delta/2}p\,$ and $q$ by $e^{\Delta_k/2}q\,$ in the above identity, it
follows that for all $p,q\in \Pi^N$
\[ c_k\int_{\b R^N} V_k p\,(x)\,q(x)\,e^{-|x|^2/2}\,w_k(x) dx\,=\, c_0
\int_{\b R^N} p(x)\, e^{-\Delta/2} e^{\Delta_k/2} q\,(x)\,
e^{-|x|^2/2} dx.\]
Corollary 3.3 now implies that
\[ \int_{\b R^N} V_k\, p(x)\,q(x) \, e^{-|x|^2/2}\,w_k(x) dx\,\geq 0
\quad \text{ for all }\, p,q\in \Pi_+^N.\]
For fixed $p\in \Pi_+^N$ we may therefore apply Lemma 3.6 with, let us
say, $\alpha = 1/4$, to the function $\,h(x):=
e^{-|x|^2/4}\,V_k\,p(x)\in C_b(\b R^N)$. This shows that 
$\,V_k\, p(x)\geq 0\,$ for all
$x\in \b R^N$ and yields the assertion.
\end{proof}

\section{Proof of the main result}

We start this section with a short discussion of the  algebras
$A_r \> (r>0)$  introduced in \eqref{(1.10)}. We should first point
out that these
are complex algebras, whereas in \cite{Du2} only series of real-valued
polynomials are considered.
It is easily checked that $A_r$  is a subalgebra of the space of
functions which are continuous on the ball $K_r$ and real analytic in
its interior: in fact, for real-valued $p\in \mathcal P_n^N$ and
$i=1,\ldots,N$ the inequality $\,\|\partial_ip\|_{\infty,
  K_1}\,\leq\,n\|p\|_{\infty, K_1}$ holds as a consequence of the Van
  der Corput-Schaake inequality, see \cite{Du2}; this allows to
  differentiate $\,f=\sum_{n=0}^\infty f_n \,\in A_r$ termwise and 
  arbitrary often.
The topology of $A_r$ is 
stronger than the topology
induced by the uniform norm on $K_r$. Notice also that  $A_r$ is not closed with respect to
$\|.\|_{\infty, K_r}$ and that $\,A_r\subseteq A_s$ with $\|.\|_{A_r}
\geq \|.\|_{A_s}$ for $s\leq r$.  The following observation is  straightforward:

\begin{lemma} $(A_r, \|.\|_{A_r})$ is a commutative Banach-$*$-algebra
   with the pointwise multiplication of
  functions, complex conjugation as involution, and with unit $1$. 
\end{lemma}

\begin{proof}  To
  show  completeness,  
let $(f^m)_{m\in \b Z_+}\,$ be a
  Cauchy sequence in $A_r$. Then for $\epsilon >0$ there exists an index
  $m(\epsilon)\in \b Z_+$ such that
\begin{equation}\label{(4.1)} 
\sum_{n=0}^\infty \|f_n^m\,-\, f_n^{m^\prime}\|_{\infty, K_r}\,<
  \epsilon\quad\text{ for }\, m, m^\prime > m(\epsilon).
\end{equation}
In particular, for each degree $n$ the homogeneous parts
$(f_n^m)_{m\in \b Z_+}$ converge uniformly on $K_r$, and hence within
  $\mathcal P_n^N$ to some $g_n\in \mathcal P_n^N$. It further follows
  from \eqref{(4.1)} that
\[ \sum_{n=0}^\infty \|g_n - f_n^m\|_{\infty, K_r} \,<\,\epsilon
\quad\text{ for }\, m>m(\epsilon).\]
Therefore $\,g:= \sum_{n=0}^\infty g_n\,$ belongs to $A_r$
with $\,\|g - f^m\|_{A_r} \to 0\,$ for  $m\to\infty.$ It is also
easily checked by a Cauchy-product argument that $A_r$ is an algebra
with $\,\|fg\|_{A_r}\,\leq\|f\|_{A_r}\cdot\|g\|_{A_r}\,$ for all
$f,g\in A_r$. The rest is clear.
\end{proof} 

We next determine the symmetric spectrum of $A_r,$  i.e. the subspace of the
spectrum $\Delta(A_r)$  given by
\[ \Delta_S(A_r):= \{\phi\in \Delta(A_r):\, \phi(\overline f)\,=\,
\overline{\phi(f)} \quad\text{ for all }\,\,f\in A_r\}.\]
As usual, $\Delta_S(A_r)$ is equipped with the Gelfand-topology. For 
$x\in K_r$ the evaluation homomorphism at $x$ is defined by
$\,\phi_x:\, A_r\to \b C,\> \phi_x(f):= f(x).$

\begin{lemma} $\, \Delta_S(A_r) = \{\phi_x : x\in K_r\,\}$, and the
  mapping $\,x\mapsto \varphi_x\,$ 
 is a homeomorphism from $K_r$  onto $\Delta_S(A_r)$. 
 \end{lemma}

\begin{proof}
It is obvious that $\phi_x$ belongs to $\Delta_S(A_r)$ for each $x\in
K_r$, with $\phi_x\not=\phi_y$ for $x\not= y$, and that the mapping
$x\mapsto\phi_x$ is continuous on $K_r$.  It remains to show that each
$\phi\in\Delta_S(A_r)$ is of the form $\phi_x$ with some $x\in
K_r$. To this end, put $\,\lambda_i:= \phi(x_i)$ for
$i=1,\ldots,N$. By symmetry of $\phi$ we have
$\lambda:=(\lambda_1,\ldots,\lambda_N)\in \b R^N.$ Moreover,
\[|\lambda|^2 \,=\, \phi(|x|^2)\,\leq\, \|\,|x|^2\|_{A_r} \,=\, r^2\,.\]
This shows that $\lambda\in K_r$. By definition of $\lambda$, the
identity $\,p(\lambda)=\phi(p)$ holds for all polynomials $p\in
\Pi^N$. The assertion now follows from the density of $\Pi^N$ in
$(A_r\,, \,\|.\|_{A_r}).$  
\end{proof}

\begin{proof}[{\bf Proof of Theorem 1.2}]   
Fix  $x\in \b R^N$ and put $r=|x|.$ Then the mapping 
\[ \Phi_x : f\mapsto V_kf(x)\]
is a bounded linear functional on $A_r$, and  Theorem 1.1 implies that
it is positive on the dense subalgebra $\Pi^N$ of $A_r$,
i.e. $\,\Phi_x(|p|^2) \geq 0$ for all $p\in \Pi^N$. Consequently, 
$\Phi_x$ is a  positive functional on the whole Banach-$*$-algebra
$A_r$. Now, by a well-known Bochner-type
representation theorem for positive functionals on commutative 
Banach-$*$-algebras (see e.g. Theorem 21.2 of \cite{FD}), there exists a
unique  measure $\nu_x\in M_b^+(\Delta_S(A_r))$ such that 
\begin{equation}\label{(4.2)}
\Phi_x(f)\,=\, \int_{\Delta_S(A_r)} \widehat
f(\phi)\,d\nu_x(\phi)\quad\>\text{ for all }\> f\in A_r,
\end{equation}
with $\widehat f$  the Gelfand transform of $f$. Denote by
$\mu_x$ the image measure of $\nu_x$ under the homeomorphism $\,
\Delta_S(A_r)\to K_r\,,\> \phi_x\to x$. Equation \eqref{(4.2)} then becomes
\[ V_k f(x)\,=\, \int_{\{|\xi|\leq|x|\}} f(\xi)\,d\mu_x(\xi) 
\quad\text{ for all}\>\, f\in A_{|x|}.\]
The normalization $V_k1=1$ implies that $\mu_x$ is a probability
measure on $\{\xi\in \b R^N: |\xi|\leq|x|\}.$ The uniqueness of
$\mu_x$ among the representing probability measures on $\b R^N$ is
clear, because identity \eqref{(1.11)} in particular determines the
(classical) Fourier-Stieltjes transform of $\mu_x$. Finally, the
transformation properties \eqref{(1.12)} follow
immediately from the homogeneity-preserving character of $V_k$ on
$\Pi^N$ and the invariance property  $V_k\circ g = g\circ V_k$  for  all $g\in G$, see Theorem 2.3 of \cite{Du2}. This
finishes the proof of Theorem 1.2. 
\end{proof}

\section{Some consequences and applications}

In this final section, we discuss only a short selection of
implications which arise from the positivity of Dunkl's intertwining
operator.  We expect that  several more useful
applications can be found, and it would of course also be of 
interest to have an explicit form for $V_k$ for larger classes of
reflection groups. In the sequel, it is always assumed that $k\geq
0$. The most prominent consequence
of Theorem 1.2, as already mentioned in the introduction, 
 is positive-definiteness of Dunkl's generalized exponential
kernel. Up to now, this has been known only in the special cases where
positivity of $V_k$ is visible from an explicit integral
representation. In particular, for the reflection group $G=\b Z_2$ on
$\b R$ and multiplicity parameter $k\geq 0$, formula \eqref{(1.9)}
shows that
\[ K_G(x,iy)\,=\, c_k \int_{-1}^1 e^{itxy}\, (1-t)^{k-1}(1+t)^kdt\,=\,
  e^{ixy}\,\phantom{}_1F_1(k,2k+1,-2ixy).\]
The following general result is an  immediate consequence of  Theorem 1.2 with $f(x) = e^{\langle x,\,z\rangle}$ and Bochner's theorem:

\begin{proposition}
For each $z\in \b C^N$, the function $\,x\mapsto K_G(x,z)\,$ has the 
Bochner-type representation
\begin{equation}\label{(5.1)}
 K_G(x,z)\,=\,\int_{\b R^N} e^{\langle\xi,\,z\rangle} d\mu_x(\xi);
\end{equation}
here the $\mu_x$ are  the representing measures from Theorem 1.2. In
particular, $K_G(x,y)>0$ for all $x,y\in \b R^N,$ and for each 
$y\in \b R^N$ the function $\,x\mapsto K_G(x,iy)\,$ is
positive-definite on $\b R^N$. 
\end{proposition}


\begin{corollary}
For each $y\in \b R^N$, the generalized Bessel function $x\mapsto J_G(x,iy)$  
is positive-definite on $\b R^N$. 
\end{corollary}

We mention that for
the group $G=S_3$ this corollary follows from the integral
representations in \cite{Du4}.
The following useful estimates of $K_G$ are immediate from 
\eqref{(5.1)}; they partially sharpen those of \cite{deJ}. (Notice that our proof of Theorem 1.2
did not involve any results from \cite{deJ}.)

\begin{corollary} Let $\nu\in \b Z_+^N$ and $|\nu| = \nu_1 +\ldots +
  \nu_N\,.$ Then for all $x\in \b R^N$ and $ z\in \b C^N,$
\begin{equation}\label{(5.4)}
|\partial_z^\nu K_G(x,z)|\le |x|^{|\nu|} \cdot e^{|x|\cdot |{\rm Re}\>
  z|}\,.
\end{equation}
In particular, $\,|K_G(x,iy)|\leq  1$ for all $x,y\in \b R^N$.
\end{corollary}

From the integral representation (5.1) we also obtain further
knowledge about the support
of the representing measures $\mu_x$:

\begin{corollary}
The measures $\mu_x\,,\> x\in \b R^N,$ satisfy \parskip=-1pt
\begin{enumerate}\itemsep=0pt
\item[\rm{(i)}]
$\,\displaystyle {\rm supp}\,\mu_x\,$ is contained in $\,{\rm co}\,\{gx,\>g\in G\},$ the convex
hull of the orbit of $x$ under $G$. 
\item[\rm{(ii)}]$\,\displaystyle {\rm supp}\, \mu_x\,\cap\,\{gx,\>g\in
  G\}\not=\emptyset. $
\end{enumerate}
\end{corollary}

\begin{proof}
(i) follows from Corollary 3.3 of \cite{deJ}. For the proof of (ii) it is
therefore enough to show that
\[ \text{supp}\,\mu_x\,\cap \,\{\xi\in \b
R^N:\,|\xi|=|x|\}\,\not=\,\emptyset.\]
Suppose in the contrary that $\,\text{supp}\,\mu_x\,\cap \,\{\xi\in \b
R^N:\,|\xi|=|x|\}\,=\,\emptyset\,$ for some $x\in \b R^N$. Then there exists a
constant $\sigma\in ]0,1[$ 
 such that
$\,\text{supp}\,\mu_x\subseteq\{\xi\in \b R^N: |\xi|\leq\sigma|x|\}.$
This leads to the estimation
\begin{equation}\label{(5.2)}
 K_G(x,y)\,=\, \int_{\{|\xi|\leq\sigma |x|\}} e^{\langle\xi,y\rangle}\,
d\mu_x(\xi)\,\leq\, e^{\sigma |x||y|}
\end{equation}
for all $y\in \b R^N$.  On the other
hand, Theorem 3.2 of \cite{Du2} with $z=0$ says that
\[ c_k \int_{\b R^N}
K_G(x,y)\,e^{-(|x|^2+|y|^2)/2}\,w_k(y)\,dy\,=\,1.\]
In view of \eqref{(5.2)}, and as $\, |w_k(y)|\leq 2^\gamma
|y|^{2\gamma}\,$ with $\, \gamma:=\sum_{\alpha\in R_+} k(\alpha)$,
it follows that 
\begin{align}\label{(5.3)} 
 1\,\leq \,&\,2^\gamma c_k\int_{\b R^N} e^{-(
  |x|^2+|y|^2)/2}\,e^{\sigma|x||y|} \, |y|^{2\gamma}dy \notag\\
  =\,& d_k \int_0^\infty e^{-(r-|x|)^2/2}\, e^{(\sigma-1)r|x|}
  r^{2\gamma+N-1}dr\notag\\
  \leq\,& d_k \int_{-\infty}^\infty e^{-r^2/2}
  e^{(\sigma-1)(|x|+|r|)|x|}\, (|x|+|r|)^{2\gamma+N-1}dr,
\end{align}
with some constant $d_k>0$ (which is independent of $x$). The
integrand in \eqref{(5.3)} is obviously majorized by $\, C\cdot
e^{-r^2/2} \,$ with a constant $C$ independent of $r$ and $x$, and
converges pointwise to $\,e^{-r^2/2}\,$ with $|x|\to \infty.$ The
dominated convergence theorem now implies that the integral
\eqref{(5.3)} tends to $0$ with $|x|\to\infty$, a contradiction. 
\end{proof} 

\noindent
We give two further applications: 

\medskip

\noindent
{\bf 1. Summability of orthogonal series in generalized harmonics.} 
The study of generalized spherical harmonics associated with a finite
reflection group and a multiplicity function $k\geq 0$ was one of the 
starting points of
Dunkl's theory in \cite{Du2} and has been extended in \cite{X2} and
\cite{X3}. Many results for classical spherical harmonics
carry over to these spherical $k$-harmonics, where harmonizity
is now meant with respect to $\Delta_k$. 
In particular, there is a natural
decomposition of $\mathcal P^N_n\vert_{S^{N-1}}$ into subspaces of
$k$-spherical harmonics, which are orthogonal 
in $L^2(S^{N-1}, w_k(x)dx).$  In \cite{X2}, Ces\`aro
summability of generalized Fourier expansions with respect to an
orthonormal basis of spherical $k$-harmonics is studied. Recall that a
sequence $\{s_n\}_{n\in \b Z_+}$ is called Ces\`aro summable of order
$\delta$ to $s$,  for short, $(C,\delta)$-summable to $s$, if 
\[ \frac{1}{\binom{n+\delta}{n}}\sum_{k=0}^n
\binom{n-k+\delta-1}{n-k}s_k\, \longrightarrow\, s \quad \text{ with
  }\, n\to\infty.\]
The following result is 
proven in \cite{X2} under the requirement that the intertwining
  operator 
$V_k$ is
positive on $\Pi^N$; Theorem 1.1 now assures its validity  for all $k\geq 0$:

\begin{theorem} Let $f: S^{N-1}\to \b C$ be continuous, and let
  $\{s_n\}$ denote the sequence of partial sums in 
    the expansion of $f$ as a Fourier series with respect to a fixed
  orthonormal basis of 
  spherical $k$-harmonics. Then $\{s_n\}$ is  uniformly
  $(C,\delta)$-summable over $S^{N-1},$ provided $\delta>\gamma+N/2-1 $ with
    $\gamma=\sum_{\alpha\in R_+} k(\alpha)$.
\end{theorem}

\medskip

\noindent
{\bf 2. Generalized moment functions.} 
Recently, in \cite{RV} a concept of Markov kernels and Markov
processes which are homogeneous
 with respect to a given Dunkl transform has been developed. In this
 context, generalized moment functions on $\b R^N$ 
 provide a useful tool. They generalize the classical monomial moment functions
 $\, m_\nu(x) = x^\nu, \>\nu\in \b Z_+^N$  and are defined as
 the unique analytic coefficients in the expansion
\[ K_G(x,y)\,=\, \sum_{\nu\in \b Z_+^N} \frac{m_{k,\,\nu}(x)}{\nu!}\, y^\nu
 \qquad(x\in \b R^N,\> y\in \b C^N).\]
From the definition of $K_G$ it follows that 
\[ m_{k,\,\nu}(x)= V_k(x^\nu) \qquad\text{for }\> \nu\in \b Z_+^N\,, \]
and Theorem 1.2 in particular implies the following useful relations
for the generalized moment functions, which are obvious only in the
classical case (again, we assume 
$k\geq 0)$:
\[ |m_{k,\,\nu}(x)|\, \leq\, |x|^{|\nu|} \quad\text{ and }\quad 0\leq
m_{k,\,\nu}(x)^2\,\leq m_{k,\,2\nu}(x) \quad\>\>\text{for all }\> x\in \b R^N,
\>\nu\in \b Z_+^N\,.\]
The first inequality is clear from the support properties of the
measures $\mu_x$ while the second one follows from Jensen's
inequality. Among the applications  of these moments, we mention
 the construction of
martingales from Dunkl-type Markov processes; for details, we refer to
\cite{RV}.

\bigskip

\noindent\emph{Acknowledgement.} It is a great pleasure to thank 
Charles F. Dunkl and Michael Voit
 for several valuable comments and  discussions. I would also like to
 express my gratitude to
  the Department of Mathematics at the University of Virgina in
 Charlottesville for its hospitality.

\end{document}